%Authors:S.J. Dilworth  and  Maria Girardi 

%Title: On Various Modes of Scalar Convergence in L_0(X) 

%Filename:dilworthgirardisclrcnvrg.atx
%TeX: AMSTeX 
%Length: 69497
%Received Date:4/5/96
%SubjectClass: 46E30
%Abstract:A sequence $\{f_n\}$ of strongly-measurable functions
% taking values in
% a Banach space $\X$ is scalarly null a\.e\.
%(resp. scalarly null in measure)
%if $x^*f_n \rightarrow0$ a\.e\. (resp. $x^*f_n \rightarrow 0$ in measure)
%for every $x^*\in \X^*$. Let $1\le p\le \infty$.
%The main questions addressed in this paper are
%whether an $L_p(\X)$-bounded  sequence that is scalarly null a\.e\.
% will converge  weakly
%a\.e\. (or have a subsequence which converges weakly a\.e\.),
% and whether an $L_p(\X)$-bounded sequence that is
%scalarly null in measure will
%have  a subsequence that is scalarly null a\.e.
%The answers to these
%and other similar questions
% depend upon $p$ and upon the geometry of $\X$. 
%

%Citation:
%

%%%%%%%%%%%%%%%%%%%%%%%%%%%%%%%%%%%%%%%%%%%%%%%%  
%%
%% 
%%         On Various Modes of 
%%     Scalar Convergence in L_0(X)     
%% 
%%      
%%    S.J. Dilworth  &  Maria Girardi
%%         
%%
%%             5  April  1996
%%
%%           AMS-TEX  version 2.1
%%
%%%%%%%%%%%%%%%%%%%%%%%%%%%%%%%%%%%%%%%%%%%%%%% 
\documentstyle{amsppt}
%\magnification=\magstep1
%\magnification=1200  %% for msri format         
%\NoRunningHeads 
\TagsOnRight
%%%%%%%%%%%%%%%%%%%%%%%%%%% MACROS   %%%%%%%%%%%%%%%%%%%%
\def\X{\frak X} 
\def\Y{\frak Y}

\def\lmr{L_0 (\Bbb R)}
\def\lpr{L_p (\Bbb R)}

\def\lmx{L_0 (\X)}

\def\lpx{L_p (\X)}
\def\lrx{L_r (\X)}

\def\e{\varepsilon}
\def\w{\omega}

\def\D{\Delta} 
\def\O{\Omega} 
\def\rn{\right\Vert}
\def\ln{\left\Vert}
\def\lav{\left\vert}
\def\rav{\right\vert}
\def\qed{\hfill{\vrule height6pt  width6pt depth0pt}}
\font\cmcsc=cmcsc8
\def\hf{\hfill}
%%%%%%%%%%%%%%%%%%%%%%%%%%%%%%%%%%%%%%%%%%%%%%%%%%%%%%%%%%%%%%%%%%%%% 
\def\ra{\longrightarrow}

\def\da{\downarrow}
\def\ane{\nearrow}

\def\ase{\searrow}

\def\nla{\longleftarrow\hskip-.9em\shortmid}
\def\nua{\uparrow\hskip-.35em\nshortmid}
\def\nda{\downarrow\hskip-.35em\nshortmid}
\def\nane{\nearrow\hskip-.58em\shortmid}
\def\nasw{\swarrow\hskip-.565em\shortmid}

\def\nanw{\nwarrow\hskip-.54em\shortmid}

\def\ranla{^{\longrightarrow}_{\nla}}
\def\span{\text{sp\,}}
\hyphenation{sca-lar-ly} 
%%%%%%%%%%%%%%%%%%%%%%%%%%%%%%%%%%%%%%%%%%%%%%%%%%%%%%%%%%%%%%%%%%%%
\topmatter
%%%%%%%%%%%%%%%%%%%%%%%%%%%
\date   
		5 April 1996                                            
\enddate
%%%%%%%%%%%%%%%%%%%%%%%%%%%%%%
\title  
On various modes of 
Scalar  Convergence in $\lmx$
\endtitle 
\author 
S. J. Dilworth 
\footnote"{$^\dagger$}"{Participant, NSF
Workshop in Linear Analysis
\& Probability, Texas A\&M University 
(supported in part by NSF grant  DMS-9311902). 
Participant,
 Convex Geometry and Geometric Functional Analysis
Program at the 
Mathematical Sciences Research Institute, 
Berkeley, CA 
(supported in part by NSF grant  DMS-9022140).
\hfill\phantom{M}}
 \quad  {\cmcsc and}   \quad
Maria Girardi $^\dagger$
\footnote"*"{Supported in part by NSF grant DMS-9306460.    
	   \hfill\phantom{P}}
\endauthor
\address \hfill\newline
Department of Mathematics, University of South Carolina, 
Columbia, SC  29208, U.S.A.  
\endaddress 
\email 
dilworth\@math.sc.edu 
and 
girardi\@math.sc.edu \endemail
\abstract  
A sequence $\{f_n\}$ of strongly-measurable functions 
 taking values in
 a Banach space $\X$ is scalarly null a\.e\.
(resp. scalarly null in measure)
if $x^*f_n \rightarrow0$ a\.e\. (resp. $x^*f_n \rightarrow 0$ in measure)
for every $x^*\in \X^*$. Let $1\le p\le \infty$. 
The main questions addressed in this paper are  
whether an $L_p(\X)$-bounded  sequence that is scalarly null a\.e\. 
 will converge  weakly 
a\.e\. (or have a subsequence which converges weakly a\.e\.),
 and whether an $L_p(\X)$-bounded sequence that is 
scalarly null in measure will
have  a subsequence that is scalarly null a\.e. 
The answers to these 
and other similar questions
often
 depend upon $p$ and upon the geometry of $\X$.  
\endabstract 
\headline={\ifodd\pageno\rightheadline\else\leftheadline\fi}
\def\rightheadline{\eightrm\hfil  Scalar Convergence \hfil\folio}
\def\leftheadline{\eightrm\hfil DILWORTH \  and \  GIRARDI \hfil\folio}
\endtopmatter

\document
\baselineskip=14pt 

\heading{1. INTRODUCTION}\endheading \vskip 10pt

Consider the space  $\lmx$ of 
strongly-measurable functions   
 defined  on  the usual % \linebreak
Lebesgue measure space $(\Omega, \Sigma, \mu)$ on $[0,1]$  
and taking  values in the Banach space~
$\X$. Among the 
most important
linear subspaces of $L_0(\X)$ are the 
Bochner-Lebesgue spaces $\lpx$ 
for $1\leq p \leq \infty$. 
When $\X = \Bbb R$,  
%and $p = 0$ or $1 \leq p \leq \infty$ (or $p=0$), 
the
$L_p(\X)$ spaces are just the usual Lebesgue spaces, 
which we shall denote by $L_p$. 
A  sequence $\{ f_n \}$  in   $\lmx$  
may converge to $f$ in $\lmx$ in 
a variety of modes.  In this paper, 
we examine  the implications 
going between the four modes described below.

The sequence $\{ f_n \}$ converges 
{\it scalarly a\.e\.\/} 
(resp.~ {\it scalarly in $L_p$\/} 
for a fixed $1 \leq p \leq \infty$, resp.~ 
{\it scalarly  in measure\/}) 
to $f$ if for each $x^*$ in the dual space ~$\X^*$ of~ $\X$,  
the corresponding sequence  $\{ x^* f_n \}$ in $\lmr$ converges  
almost everywhere (resp.~  in  $\lpr$, resp.~ in measure)  to $x^*f$.  
Thus, $\{ f_n \}$ converges scalarly a\.e\.  to $f$ 
if for each $x^*$ in $\X^*$ there is a   
set $A$ (which depends on $x^*$)
of full measure  such that
$\{ x^* f_n (\omega) \}$ converges to 
$x^* f(\omega)$ for each $\omega$ in $ A$. 
If the sequence satisfies the stronger 
property that there is a set $A$ (independent of  $x^*$) 
of full measure  such that
$\{ x^* f_n (\omega) \}$ converges to 
$x^* f(\omega)$ for each $\omega$ in $A$ and each 
$x^*$ in $\X^*$, then we say that 
$\{ f_n \}$ converges 
{\it weakly almost everywhere\/} 
(i\.e\.~{\it weakly a\.e\.\/}) to $f$. 

The following obvious positive and negative implications hold:    
$$\vbox{
\settabs\+\indent
&\hskip 1 true in &\hskip 2 true in &\hskip 1 true in& \cr 
\+&\hf weakly a\.e\.\hf&\hf$\ra$\hf&\hf scalarly a\.e\. \hf\cr
\+\cr
\+&\hf$\nua~\nda$\hf&\hf$\nane~\nasw\quad\ase~\nanw$\hf&\hf$\nua~\da$\hf&\cr
\+\cr
\+&\hf scalarly in $L_p$ \hf&\hf$\ranla$\hf&\hf scalarly in measure \cr}$$ 
for each $1\leq p < \infty$. Similarly:    
$$\vbox{
\settabs\+\indent
&\hskip 1 true in &\hskip 2 true in &\hskip 1 true in& \cr 
\+&\hf weakly a\.e\.\hf&\hf$\ra$\hf&\hf scalarly a\.e\. \hf\cr
\+\cr
\+&\hf$~\nda$\hf&\hf$\ane~\nasw\quad\ase~\nanw$\hf&\hf$\nua~\da$\hf&\cr
\+\cr
\+&\hf scalarly in $L_p$ \hf&\hf$\ranla$\hf&\hf scalarly in measure \cr}$$ 
for $p = \infty$.

Sets of the form
$$
V_{x^*, \e}(0) = \left\{ g\in L_0(\X) \: 
  \mu \left\{ \lav x^* (g) \rav \geq \e \right\} < \e \right\} \ ,    
$$ 
where $x^* \in \X^*$ and $\e >0$,  
form a local 
subbasis at zero for the  
translation-invariant topology of 
scalar convergence in measure. 
Endowed with this topology,  $L_0(\X)$  
is a non-locally convex  Hausdorff 
topological vector space. 
For a fixed $1\leq p \leq \infty$, 
sets of the form 
$$ 
S_{x^* , \e}(0) = \left\{ g \in L_p (\X) \:
     \ln x^* g \rn_{L_p}  < \e \right\} \ , 
$$ 
where  $x^* \in \X^*$ and $\e>0$, form a local 
subbasis at zero for the  
translation-invariant topology of 
scalar convergence in $L_p$.  
Endowed with this topology, $L_p(\X)$ 
is a locally convex  Hausdorff 
topological vector space.

Let $\X_0$ be a norm-closed subspace of $\X$.
The Hahn-Banach Theorem quickly gives two observations. 
First,
a sequence of  $L_0(\X_0)$ functions 
converges to the null function in one of the above modes when viewed as 
a sequence in $L_0(\X_0)$  if and only if it does   
so when  viewed as a sequence in $L_0(\X)$. 
Secondly,  the topology of scalar convergence in measure
on $L_0(\X_0)$  coincides with the subspace topology inherited from  
the topology of scalar convergence in measure on
$L_0(\X)$. Let us show that under this identification 
 $L_0(\X_0)$ is in fact a {\it closed}
subspace of $L_0(\X)$.
\proclaim{Proposition 1.1} 
Let $\X_0$ be a (norm-closed) subspace of $\X$. Then
$L_0(\X_0)$ is a closed subspace of $L_0(\X)$
in the topology of scalar convergence in measure.
In particular, if
$\{f_n\}$ is a sequence of $\X_0$-valued functions in $L_0(\X)$  
that converges 
scalarly in measure to $f$ in  $ L_0(\X)$, then $f$ is  also $\X_0$-valued. 
\endproclaim
\demo{Proof} Let $f\in L_0(\X)$ belong to the closure of $L_0(\X_0)$ 
in the topology of scalar convergence in measure. 
Since the range of $f$ is essentially separably-valued,  there
is a subset $Y \supset \X_0$ such that $Y/\X_0$ is separable and
$f(\omega) \in Y$ a\.e\. By the Hahn-Banach theorem there exists
a sequence $\{x_n^*\}$ in $\X^*$ such that if $y \in Y$ then $y \in \X_0$ 
if and only if $x_n^*(y)=0$ for all $n \ge 1$. In particular, if
$g \in L_0(\X_0)$ then $x_n^*(g(\omega))=0$  for all $n \ge1$ a\.e\.,
and it follows that $x_n^*(f(\omega))=0$  for all $n \ge 1$ a\.e\.,
which proves that $f \in L_0(\X_0)$. \qed \enddemo

For $1 \leq p \leq \infty$,         
it is easy to see that 
$L_p (\X)$ is {\it not} a closed  
subspace of $L_0 (\X)$ in the topology of  scalar convergence in measure 
for any Banach space $\X$.
So we consider the unit ball 
$$B(L_p(\X)) = \{ g\in L_0(\X) \: \ln g \rn_{L_p(\X)} \leq 1 \}$$ 
of $L_p(\X)$.  
   
\proclaim{Proposition 1.2} For $1\le p\le \infty$,   
the  unit ball $B(L_p(\X))$ of $L_p(\X)$ is
closed in $L_0(\X)$ in the topology of scalar
convergence in measure. 
\endproclaim

\demo{Proof}  
Consider $f \in L_0(\X)\setminus B(L_p(\X))$.    
It is sufficient to find an open neighborhood about $f$ 
that does not meet $ B(L_p(\X))$.
   
Fix $\e > 0$ so that $(1-4\e) \ln f \rn_{L_p(\X)} > 1$. 
Since $f$ is strongly-measurable, 
there is a countably-valued function 
$g\in L_0(\X)$ satisfying 
$$\ln f(\omega) - g(\omega) \rn_{\X}~\leq~ 
\e~ \max\{\ln f(\omega)\rn_{\X},\ln g(\omega)\rn_{\X}\}$$   
for almost all $\omega$.  
By making an appropriate choice of representative we may
write $g=\sum_k x_k 1_{E_k}$, where $\{E_k\}_k$ 
partitions the support of $g$  
into sets of strictly positive
measure. Now
$\ln g \rn_{L_p(\X)} \ge (1-\e) \ln f \rn_{L_p(\X)}$, and so
$(1 -3\e)\ln g \rn_{L_p(\X)}>1$. Hence we may choose $N \in \Bbb N$
so that
$$(1 -3\e)\ln \tilde g \rn_{L_p(\X)}>1~,$$
where $\tilde g= \sum_{k=1}^N x_k 1_{E_k}$.
Now find  $\{ x^*_k \}_{k=1}^N$ in $S(\X^*)$ 
with  $x^*_k (x_k) = \ln x_k \rn $. 

Consider the following neighborhood of $f$ in the topology
of scalar convergence in measure:
$$U=
 \bigcap_{k=1}^N \left\{ h\in L_0(\X) \: 
 \mu \left\{ \lav x^*_k (h-f) \rav \geq \e \ln x_k \rn \right\} 
   < \e \mu(E_k) \right\}~.  
$$
For $1 \leq k \leq N$, 
if $\omega \in E_k$   
then % $\tilde g(\omega) = \|x_k\|$ and also 
$\|f(\omega) - \tilde g(\omega)\| \le \e \| x_k \|$; 
thus,  for 
each $h \in U$ 
$$
\mu \{ \omega \in E_k: 
  \|h(\omega)\| \ge (1 - 2\e) \| x_k \| \}
  \ge (1-\e)\mu(E_k)~.
$$
Hence 
$$
\ln h \rn_{L_p(\X)} \ge (1-\e)^{1/p}(1-2\e)
\ln \tilde g \rn_{L_p(\X)}>(1-3\e)\ln \tilde g \rn_{L_p(\X)}>1~,
$$
following the convention that $1/\infty = 0$.  
So $U$ does not intersect $ B(L_p(\X))$, as required.
\qed
\enddemo

This suggests imposing the following 
natural boundedness  conditions:   
a sequence $\{ f_n \}$  of $L_p(\X)$ functions is said to be: 
\roster
\item"{-}"  {\it pointwise bounded a\.e\.}  if 
$\sup_n \ln f_n (\omega) \rn_{\X}<\infty$  
for each $\omega$ in some set of full measure.     
\item"{-}"  
 bounded  in $L_p(\X)$ 
(for short, {\it $L_p(\X)$-bounded}) if  
 $\sup_n \ln f_n \rn_{L_p(\X)}<\infty$.  
\endroster  From Proposition 1.2, we see that
sequences of the latter type are well-behaved in the following sense.  

\proclaim{Corollary 1.3} Let $1\le p\le\infty$ and 
suppose that $\{f_n\}$ is bounded in $L_p(\X)$.
If $\{f_n\}$ converges scalarly in measure to $f\in L_0(\X)$, then
$f\in L_p(\X)$. \endproclaim
\demo{Remarks}    
1. Hence, in discussing the convergence (in any one of the above four modes) 
of a sequence in $L_0(\X)$ of functions valued in a subspace $\X_0$ of $\X$, 
there is no loss of generality in taking the limit function to be the null function   
and viewing the sequence as  in $L_0(\X_0)$.  
\newline  
2. Similarly, if we choose to restrict ourselves 
to the subset $L_p(\X)$ of $L_0(\X)$, 
in discussing scalar convergence in measure
for an $L_p(\X)$-bounded sequence,    
   there will be no loss of 
generality in taking the limit function to be the null function.
\newline
3. The question of the {\it existence} of a limit for an $L_p(\X)$-bounded
{\it Cauchy sequence} in the topology of scalar convergence of measure
is more problematic and will be deferred until Section~5.
\enddemo

We will also use the following  elementary facts without further comment. 
Fact~1.4 provides a useful necessary condition for weak a\.e\. 
convergence while Fact~1.5 will be used to prove scalar  
convergence a\.e\.
\demo{Fact 1.4} A weakly convergent sequence in a Banach space 
is norm-bounded.
Thus, %Hence,
if for a given sequence $\{ f_n \}$ in $L_0(\X)$, 
there exists a subset $B$ of strictly positive 
$\mu$-measure such that $\limsup \ln f_n (\omega) \rn = \infty$ 
for each $\omega \in B$, 
then   $\{ f_n \}$ does not converge weakly a\.e\. 
\enddemo
\demo{Fact 1.5}
A sequence $\{ f_n \}$ in $L_1$  
converges to the null function  
a.e. whenever 
$\sum \ln f_n \rn_{L_1} < \infty$.  
\enddemo

If $Y$ is a subset of $\X$, 
then $\span  Y $ denotes the linear span of $Y$ 
and  $\left[ Y \right]$ denotes 
the closed linear span of $Y$.    
All notation and terminology, not otherwise explained,  
are as in [DU] or [LiT].

\heading{2. CONVERGENCE PROPERTIES}\endheading \vskip 10pt

Proposition~1.2 suggests that it would be of
interest to study the
following properties 
that a Banach space $\X$ might enjoy.   
\roster
\item"($A_p^\prime$)" Each  $L_p(\X)$-bounded  
sequence of functions 
that converges  scalarly a\.e\. 
 also converges weakly a.e.
\item"($A_p$)" Each  $L_p(\X)$-bounded  
sequence of functions 
that converges  scalarly a\.e\. 
 has a subsequence 
 that  converges weakly a.e.
\item "($B_p$)" Each  $L_p(\X)$-bounded    
sequence of  functions 
that converges  scalarly in measure 
has a subsequence that converges weakly a.e. 
\item"($C_p$)" Each  $L_p(\X)$-bounded  
sequence of functions 
that converges   scalarly in $L_p$  
has a subsequence that converges weakly a.e.
\item"($D_p$)" Each  $L_p(\X)$-bounded  
sequence of functions  
that converges  scalarly in measure  
has a subsequence that converges 
scalarly a.e. 
\item"($E_p)$" Each $L_p(\X)$-bounded sequence
of functions that converges scalarly in $L_p$
has a subsequence that converges scalarly a.e.
\endroster
For convenience,  a 
schematic  summary of the properties is given below, in which 
a double arrow indicates  an implication 
that is always valid. 
$$
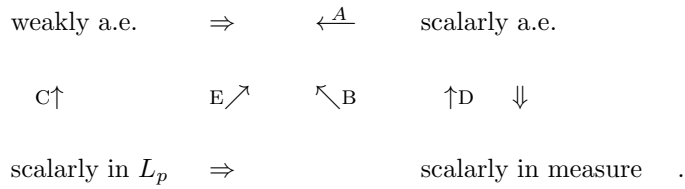
\vbox{
\settabs\+indent
&\hskip 75 pt
&\hskip 40 pt
&\hskip 40 pt
&\hskip 60 pt\cr 
\+&weakly a\.e\.
  &$\Rightarrow$ 
  &$\longleftarrow^{\hskip-1 em A}$ 
  &scalarly a\.e\.\cr
\+& \cr
\+&\hskip 9 pt{\smc c}$\uparrow$
  &{\smc e}$\nearrow$
  &$\nwarrow${\smc b}
  &\hskip 9 pt $\uparrow${\smc d}
   \hskip 11 pt $\Downarrow$  \cr
\+& \cr    
\+&scalarly in $L_p$
  &$\Rightarrow$ 
  &
  &scalarly in measure \quad . \cr
}
$$

\demo{Remarks} 1.   
Clearly, subsequential convergence 
is the most one can expect 
in passing from scalar in measure   or scalar in $L_p$
 to scalar a.e. convergence.\newline
2. 
Note that if $\X$ has  (Property$_p$) 
and $p < q$, then $\X$ also has (Property$_q$). \newline 
3. Note the following obvious implications.
$$\vbox{
\settabs\+\indent&\hskip 30 pt &\hskip 30 pt &\hskip 30 pt &\hskip 30 pt &\hskip 30 pt\cr
\+&($A_p$)&$\longleftarrow$&($B_p$)&$\longrightarrow$&($C_p$)\cr
\+&       &              &\hskip 7 pt$\downarrow$ &           &\hskip 7 pt$\downarrow$\cr
\+&       &              &($D_p$)&$\longrightarrow$&($E_p$)\cr }
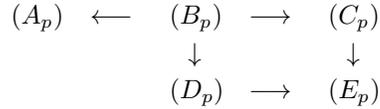
$$ 
\newline  
4. For a fixed $1 \leq p < \infty$, a  sequence $\{f_n\}$ in $L_0(\X)$ 
converges scalarly in $L_p$   
if and only if 
(i) $\{f_n\}$ converges scalarly in measure and 
(ii) for each $x^* \in \X^*$, the set 
$\{ \lav x^* f \rav^p \}$ is uniformly integrable.  
Note that (ii) holds when $\{f_n\}$ is $L_r(\X)$-bounded 
for some $r > p$.     
\enddemo

\heading{3. $L_\infty(\X)$-BOUNDED  SEQUENCES}\endheading \vskip 10pt
\subhead [3.i] From scalar convergence to weak a.e. convergence 
\endsubhead 
In this subsection, we characterize those spaces having 
($A_\infty$),  ($B_\infty$), and ($C_\infty$).

\proclaim{Theorem~3.1} 
Let $\X^*$ have the 
Radon-Nikod\'ym property. Then 
$\X$ enjoys  the following properties: 
\roster 
\item   
 Let $\{f_n\}$ be a sequence in $L_0(\X)$ which converges  scalarly a\.e\.
 to some 
 $f \in L_0(\X)$.
 Then $\{f_n\}$ converges weakly a\.e\. to $f$ 
 if and
 only if $\{f_n\}$ is pointwise-bounded a\.e\.
 \item
 Let $\{f_n\}$ be a sequence in $L_0(\X)$ which converges
 scalarly a\.e\. to some $f \in L_0(\X)$. 
 Then $\{f_n\}$ has a subsequence which 
 converges weakly a\.e\. to $f$ if and
 only if $\{f_n\}$ has a subsequence which    
 is pointwise-bounded a\.e\.
 \item
 Let $\{f_n\}$ be a sequence in $L_0(\X)$ which converges
 scalarly in measure. Then $\{f_n\}$ has a subsequence which 
 converges weakly a\.e\. if and
 only if $\{f_n\}$ has a subsequence which  
 is pointwise-bounded a\.e\.
\endroster   
\endproclaim  
   
\demo{Proof}   
By Fact 1\.4 each a\.e\. weakly convergent sequence is bounded a\.e\.,
and so necessity in (1)-(3) (for an {\it arbitrary} Banach space)
is clear. To prove sufficiency  observe that
we may take $f=0$ without loss of generality. We may also assume that
$\X^*$ is separable.
(Indeed, by the Pettis measurability   
theorem [P], there is a separable subspace $\X_0$ of $\X$ such that 
the $f_n$'s are essentially valued in $\X_0$.  Because  $\X^*$ has the RNP, 
$\X^*_0$ must be separable [S, Thm.~2].)      
Let $\{ x_i^* \}$ be dense in $\X^*$. Now we prove sufficiency 
in (1).
Since $\{ f_n \}$ is scalarly null a\.e\., it follows that 
 for each $i$ there is a set $A_i$ of full 
measure such that if $\omega\in A_i$, then 
$\lim_n x^*_i f_n(\omega) = 0$.  Put $A = \cap_i A_i$.  Since 
the $f_n$'s are pointwise-bounded  on some set $B$ of full measure, 
and since 
$\{ x_i^* \}$ is dense in 
$\X^*$, it follows that   
$\lim _n x^* f_n (\omega) = 0 $ 
for each $x^* \in \X^*$ and for each $\omega\in A\cap B$.  
Thus, $\{ f_n \}$ is 
weakly  null a\.e\. Sufficiency in (2) follows at once.  
Finally, we prove sufficiency in (3). By first  
passing to a subsequence we may assume that $\{f_n\}$
is pointwise bounded almost everywhere.  
Since $\{ f_n \}$ converges   to zero scalarly in measure,    
 for each $i$  
 the sequence 
$\{ x^*_i f_{n} \}_n$ converges in measure to the null 
function.  So by a Cantor diagonalization argument there exists 
a subsequence
$\{ f_{n_k} \}$ such that for almost all $\omega$ 
$$ 
  \lim_k x^*_i \left( f_{n_{k}} \left(\omega\right)\right) = 0 \ . 
$$  
for all $i$. Now, arguing as before, the pointwise boundedness 
implies that $\{ f_{n_k} \}$ is weakly null a\.e\. 
\qed
\enddemo  

For the Banach spaces  $\ell_1$, $C(\Delta)$, and the James tree space,
Davis and Johnson [DJ] constructed examples 
of $L_\infty(\X)$-bounded sequences  that converge scalarly a\.e\. 
but not weakly a.e.  
They conjectured that such a sequence exists for any space $\X$ whose dual 
fails the Radon-Nikod\'ym property (RNP).  
Combined with work of Uhl [U] and Stegall [S],  a result of Edgar [E]
shows that their conjecture was correct.   
In fact, rather more can be said as the following theorem (whose 
proof was inspired by [E]) shows.

\proclaim{Theorem 3.2}  For a Banach space $\X$, the following are equivalent:  
\roster
\item  $\X^*$ has the  Radon-Nikod\'ym property (i\.e\. $\X$ is an Asplund space);
\item $\X$ has ($A_\infty^\prime$); 
\item $\X$ has ($A_\infty$); 
\item $\X$ has ($B_\infty$).
\endroster                                    
\endproclaim    
\flushpar

\demo{Proof} 
Several 
 implications  follow from 
Theorem~3.1.   
To prove the other 
implications, suppose that  
 $\X^*$ fails the RNP. 
Then [U]  
there is a separable subspace $\X_0$ of $\X$  
such that $\X_0^*$ is not separable.     
We shall construct an  $L_\infty(\X)$-bounded  sequence $\{ g_n \}$   
of $\X_0$-valued  functions such  that 
$g_n \to 0$  scalarly a\.e\. and scalarly in $\lrx$ for  
$1 \leq r < \infty$,   
but such that no subsequence of  $\{g_n\}$ converges weakly a\.e.
(This particular construction has  been fruitful  
in several similar characterizations of 
$\X^*$ having the RNP [e\.g\.~E, GS, DG].)  
This will show that $\X$ fails ($A_\infty$),
thus completing the proof of the theorem.

Let $\Delta = \{ -1, 1 \}^{\Bbb N}$  be the Cantor group with  Haar measure $\nu$. 
Let $\{ \Delta^n_k \:  k=1,\ldots, 2^n \}$  
be the standard $n$-th partition of $\D$. 
Thus  $\D^0_1 = \D$ and $\D^n_k = \D^{n+1}_{2k-1} \cup \D^{n+1}_{2k}$ 
and $\nu(\D^n_k) = 2^{-n}$.
Instead of our usual Lebesgue measure space on $[0,1]$,    
we shall now take our underlying  measure space to be the completion of
$\nu$ for  the completion of the 
Borel $\sigma$-algebra of $\Delta$.
 Thus
$L_p(\X)$ will denote $L_p(\D,\nu;\X)$.

We
consider the space $C(\D)$  of  real-valued   
continuous functions on $\D$ as a subspace of $L_\infty (\D,\nu)$. 
Let 
$ \{ 1_{\D} \}  \cup 
  \{ h^n_k \:  n=0,1,2, \ldots \text{  and  } k=1,\ldots, 2^n \}$  
be the usual Haar basis of $C(\D)$, where $h^n_k \: \D \to \Bbb R$ is given by 
$$
    h^n_k ~=~ 1_{\D^{n+1}_{2k-1}} ~-~ 1_{\D^{n+1}_{2k}} \ .
$$  
Let $\{ e^n_k \:  n=0,1,2, \ldots \text{  and  } k=1,\ldots, 2^n \}$  
be  an  enumeration  (lexicographically) of the usual $\ell_1$ basis   
and let $H \: \ell_1 \to L_\infty$ be the Haar operator that takes 
$e^n_k$ to $h^n_k$. 

By Stegall's Factorization Theorem [S, Theorem~4], 
$H$ factors through $\X_0$, i\.e\. 
there are bounded linear operators 
$R \: \ell_1 \to \X_0 $ and $S \: \X_0 \to L_\infty$ 
such that $ H = SR$.  
$$\vbox{
\settabs\+\indent 
&\hskip 1 true in &\hskip 1 true in &\hskip 1 true in& \cr
\+&\hf$\ell_1\quad$
  &\hf$\longrightarrow^{^{\hskip-1.25em H}}\hskip1em$\hf
  &$L_\infty(\D)\supset C(\D)$\hf&\cr 
\+\cr
\+&\hf$\searrow_{^{\hskip-1.3em R}}$& 
     &$\nearrow_{^{\hskip-.1em S}}$\hf&\cr  
\+\cr
\+& &\hf$\X_0$\hf&\cr
}$$ 
Let $\tilde R$ be the natural extension of $R$ to 
a bounded linear operator from $L_1(\ell_1)$ to $L_1(\X_0)$. 

Consider the sequence $\{ f_m \}$ of $L_1(\ell_1)$ functions 
given by 
$$
   f_m  (\cdot) ~=~ 
	   \frac{1}{m}  \sum_{n=1}^m \sum_{k=1}^{2^n} h^n_k(\cdot) e^n_k \ .
$$
Let $g_m = \tilde R (f_m)$.  
Clearly, $\{ g_m \}$ is $L_\infty(\X)$-bounded since   
$|| f_m ( \w ) ||_{\ell_1} = 1$ for $\nu$-a.e. $\w$ and each $m$. 

To examine the scalar behavior of  $\{ g_m \}$, 
 note that if   
$y^*\in \X_0^*$, then   
$$y^* g_m(\cdot) = (R^* y^*) f_m (\cdot)~,$$ 
where $R^* y^* \in \ell^*_1$.   
So to show that 
 $\{ g_m \}$ converges to the null function 
scalarly a\.e\. and scalarly in $L_r (\X)$ 
we need only show the same for $\{ f_m \}$. 
So  fix a functional $x^*$ in $\ell^*_1$; let $x^*$ have the form 
$(\alpha^n_k) \in \ell_\infty$, lexicographically ordered.  
Then 
$$
x^* f_m(\omega) = \frac{1}{m} \sum_{n=1}^m X_n(\cdot) 
\text{\qquad where\qquad} 
X_n(\cdot) = \sum_{k=1}^{2^n} h^n_k(\cdot) \alpha^n_k \ . 
$$  
Note that $\|X_n\|_\infty \le \|x^*\|$, that 
each $X_n$ has zero mean, and that $\int X_n X_m \, d\nu = 0$ when $n\neq m$.  
The Strong Law of Large Numbers for uncorrelated random variables 
with uniformly bounded second moments [cf.~C, Thm.~5.1.2] 
 gives 
that $\{ x^* f_m \}$ converges to the null function  a\.e. 
Since $\|x^*f_m\|_\infty\le\|x^*\|$ it also follows that $\{x^*f_m\}$
converges to the null function in $L_p$ for $0\le p<\infty$.

We shall now show that no subsequence of $\{g_m\}$ converges
weakly a\.e\. Since $\{g_m\}$ is scalarly null a\.e\., 
it suffices to show that no subsequence
is weakly null a\.e\.
For $\omega\in\Omega$,  let 
$\ell_\omega\in\left[ C\left(\Delta\right)\right]^*$ 
be the point evaluation at $\omega$ functional 
and let    
$\widetilde \ell_\omega\in\left[ L_\infty\left(\Delta\right)\right]^*$  
be any Hahn-Banach extension.  
Then    
$$
\left(S^*\widetilde\ell_\omega\right) \left( g_n \left(\omega\right)\right) 
=
\ell_\omega \left( H f_n \left(\omega\right)\right)= \left(
\frac{1}{m} \sum_{n=1}^m \sum_{k=1}^{2^n} h^n_k(\omega) h^n_k(\cdot)\right)
(\omega) =1~,  $$  
and thus  no 
subsequence of  $\{g_n\}$  converges weakly a\.e\.  (to the null function)
in $L_0(\X_0)$. 
\qed
\enddemo

Property $C_\infty$, on the other hand, is a much weaker property 
according to the 
 following mildly surprising result.  
 
\proclaim{Proposition 3.3} 
Let $\{f_n\}$ be a sequence in $L_0(\X)$
that converges to the null function scalarly in $L_\infty$.
Then $\{f_n\}$ is weakly null a.e. In particular, every
Banach space enjoys ($C_\infty$)
 (and {\it a fortiori}   ($E_\infty$)).
\endproclaim 

\demo{Proof} 
For each $n\ge1$, we may
 write $f_n = g_n + h_n$,  
where $h_n$ has $L_\infty (\X)$-norm at most $1/n$ 
and $g_n$ is countably-valued  (see e\.g\. [~DU II.1.3]). 
By choosing a suitable representative of $g_n$ in $L_\infty(\X)$,
we may express
$g_n$ as 
$$ g_n = \sum_{k = 1}^\infty  x^n_k ~1_{E^n_k}, 
$$ 
where each $E^n_k$ has strictly positive measure 
and, for each $n$, $\Omega$ is the disjoint union of $\{ E^n_k \}_k$. 
Note that $\{ g_n \}$ also  
converges to the null function scalarly in $L_\infty (\X)$. Hence, for 
each $\omega\in\Omega$ and each
$x^*\in X^*$, we have
$$|x^*(g_n(\omega))|\le \sup_k |x^*(x^n_k)|=\|x^*(g_n)\|_\infty
\rightarrow 0$$
as $n\rightarrow \infty$. Hence $\{g_n(\omega)\}$ is weakly null
for all $\omega\in\Omega$. Clearly, $\{h_n(\omega)\}$ is norm-null
a.e., whence $\{f_n\}$ is weakly null a.e. \qed
\enddemo
Perhaps the following theorem is the most useful analogue
for scalar convergence in general 
Banach spaces of the familiar fact from real analysis that
each sequence  that  converges in measure has a subsequence that 
converges a.e.
\proclaim{Theorem 3.4} Let $K$ be a weakly compact subset   
of a Banach space $\X$ and let $\{f_n\}$ be
 a   sequence in $L_0(\X)$ such that each $f_n$ is 
essentially $K$-valued.
If  $\{f_n\}$
converges scalarly in measure to $f\in L_0(\X)$,
then some subsequence $\{f_{n_k}\}$ converges weakly a\.e\.
to $f$. (In particular, $f$ is essentially $K$-valued.)      
\endproclaim
\demo{Proof} Clearly, we may assume that $\X$ is separable. 
First, we show that $f(\omega)\in K'=\overline{\text{conv}(K)}$ a.e. 
Suppose this is not the case.
 Then there exists a closed ball $B\subset \X\setminus
K'$ such that $\mu(A)>0$, where $A=\{\omega:f(\omega)\in B\}>0$. By the
Hahn-Banach Theorem there exists $x^*\in\X^*$ and $\alpha,\beta
\in\Bbb R$
such that 
$$\sup_{k\in K'}x^*(k)<\alpha<\beta<\inf_{b\in B}x^*(b)~.$$
Since each  $f_n$ is (without loss of generality) 
$K$-valued, it follows that
$$ x^*(f_n(\omega))<\alpha<\beta<x^*(f(\omega))$$   
for all $n\ge1$ and all  $\omega\in A$. This contradicts
the fact that $\{f_n\}$ converges scalarly in measure to $f$.
Hence, by replacing $f_n$ by $f_n-f$ and $K$ by $K'-K'$,
we may assume without loss of generality that $\{f_n\}$
converges scalarly in measure to the null function and that
$K$ is a separable weakly compact set containing zero. 
It is easily seen that the weak topology on $K$ is
generated by a sequence $\{x_n^*\}$ in $\X^*$. By a 
Cantor diagonal argument there exists a subsequence $\{f_{n_k}\}$
and a set $\Omega'\subset \Omega$ of full measure such that 
$x^*_n(f_{n_k}(\omega))\rightarrow0$ as $k\rightarrow\infty$
for all $n\ge1$ and for all $\omega \in \Omega'$  
Now let $x^*\in \X^*$ and let $\varepsilon>0$. There exists
$\delta > 0$ and $N\ge1$ such that
$$\{k\in K:|x^*(k)|<\varepsilon\}\supset
 \cap_{i=1}^N\{k\in K:|x_i^*(k)|<\delta\}~.$$
It follows that 
$$ \{\omega\in \Omega:|x^*(f_{n_k}(\omega))|<\varepsilon\quad
\text{for all sufficiently large k}\}\supset
\Omega'~,$$
 and so $x^*(f_{n_k}(\omega))\rightarrow 0$
as $n\rightarrow \infty$ for all $\omega\in \Omega'$. Hence
$\{f_{n_k}\}$ is weakly null a.e. \qed
\enddemo
Minor variations in the above proof gives the following result.     
\proclaim{Theorem 3.5} Let $K$ be a weakly compact subset 
of a Banach space $\X$  and let $\{f_n\}$ be
 a sequence in $L_0(\X)$ 
 such that each $f_n$ is 
essentially $K$-valued.     
If $\{f_n\}$
converges scalarly a.e. to $f\in L_0(\X)$,
then  $\{f_n\}$ converges weakly a.e.
to $f$.  \endproclaim

\subhead [3.ii] From scalar convergence in measure    
to scalar a.e. convergence
\endsubhead 
In the previous subsection, 
Proposition~3.3 shows that each 
Banach space enjoys ($E_\infty$).  
In this subsection, we explore ($D_\infty$).

For each $1 \leq p \leq \infty$,
it follows directly from the definitions 
that
  ($B_p$) implies ($D_p$); 
however, they are not equivalent. 
Indeed, Theorem~3.2 implies that $\ell_1$ fails 
($B_p$) for each  $1 \leq p \leq \infty$. 
However, $\ell_1$ has ($ D_p$) for each $1 \leq p \leq \infty$, 
  as the 
next theorem, which was pointed out to us by W.B.~Johnson 
[cf.~WBJ],  shows.  
We  are grateful   to him for permission to include this 
result here.

\proclaim{Theorem 3.6} 
A scalarly null in measure sequence of $L_0(\ell_1)$ functions 
contains a scalarly null a.e.~subsequence.  
 In particular, 
$\ell_1$ satisfies properties ($D_p$) for each 
$1\leq p \leq \infty$.
\endproclaim 

The following lemma will be used in the proofs of several  
results, including Theorem 3.6.

\proclaim{Lemma 3.7} 
Let $\X$ be a Banach space with a basis  and 
 $\{ f_n \}$ be a scalarly in measure null sequence 
of $\lmx$ functions.   
There exists a subsequence $\{ f_{n_k} \}$ of 
 $\{ f_n \}$, a blocking $\{ \X_k \}$ of the basis, 
and sequences  $\{ g_k\}$ and $\{h_k\}$ 
of $\lmx$ functions  so that:  
\roster
\item"{(i)}" $\{h_k\}$  converges a.e (in $\X$-norm) to the null function; 
\endroster 
and for each $k$: 
\roster 
\item"{(ii)}"  $f_{n_k}  = g_k + h_k$;  
\item"{(iii)}" $g_k = P_k f_{n_k},$   
 where $P_k$ is the natural projection  
  of $\X$ onto $\X_k$. 
\endroster
In particular, $\{ g_k \}$ is also scalarly null in measure. 
\endproclaim 

\demo{Proof of Lemma 3.7} 
Let $\{ x_n \}_{n\geq 1}$ be a normalized basis for $\X$ 
and let $\{ x^*_n \}$ be the corresponding biorthogonal functionals.    
Consider a sequence $\{ f_n \}_{n\geq 1}$  of $L_0(\X)$ 
functions that is scalarly null in measure.
It suffices to  construct inductively two increasing sequences
 $\{ n_k \}_{k\ge1}$ and $\{ m_k \}_{k\ge0}$
of integers and a sequence  $\{ g_k \}$ of functions 
such that, 
for $\X_k \equiv \span\{x_i: m_{k-1} < i \le m_k \} $, each 
 $g_k$  satisfies (iii) and 
$$\mu(\{\omega: \|g_k(\omega)-f_{n_k}(\omega)\|\ge 2^{-k}\})\le 2^{-k}. \tag1$$
To start the induction set $m_0=0$ 
and $n_0 = 0$.   Suppose that $k\ge1$ and that
$n_i$ and  $m_i$ have been chosen  
 for $i\le k-1$. 
Since $\{ f_n \}$ is assumed to be scalarly null in measure, 
it follows that, for
each fixed $i$, the sequence $\{ x_i^*(f_n) \}_n$ 
converges to zero in measure. So there exists $n_k > n_{k-1}$
such that 
$$
\mu\left( \left\{ 
  \omega: \sum_{i=1}^{m_{k-1}}|x_i^*(f_{n_k})|\ge 2^{-k-1}
\right\} \right) 
\le 2^{-k-1} \  .
$$
Hence there exists $m_k > m_{k-1}$   such that, 
for $\X_k \equiv \span\{x_i: m_{k-1} < i \le m_k \} $,   
the function $g_k$ as given in (iii) satisfies  (1), which completes the induction.
\qed
\enddemo

\demo{Proof of Theorem 3.6} 
Consider a sequence $\{ f_n \}$  of $L_0(\ell_1)$ 
functions that is scalarly null in measure.
Let $\{ e_n \}$ be the standard basis of $\ell_1$. 
Find a blocking $\{ \X_k\}$ of $\{ e_n \}$ 
and  sequences $\{ g_k \}$ and $\{h_k \}$  as given by Lemma~3.7. 

In view of (i), it is  enough to show that
$\{ g_k \}$ has a subsequence that is scalarly null 
a\.e. 
With that in mind, we establish the following claim.
\newline
{\bf Claim.} Given $\varepsilon>0$,
$$
    \underset \|x^*\|\le1 \to\sup
    \mu\left(\left\{ |x^*(g_k)|>\varepsilon\right\}\right)\rightarrow0
$$
as $k\rightarrow\infty$. \newline
{\it Proof of Claim.} Suppose not. Then there exist
 $\varepsilon>0$ and $x_k^*
\in \ell_1^*$,
with $\|x_k^*\|\le1$, such that
$$ \limsup_k \mu(\{|x_k^*(g_k)|>\varepsilon\})>\varepsilon.$$
Note that  $\ell_1=(\sum\oplus E_k)_1$
and so $\ell_1^*=(\sum\oplus E_k^*)_\infty$. Thus there exists 
$x^* \in \ell_1^*$ such that, for each $k$, the  functionals 
 $x_k^*$ and $x^*$ 
have identical restrictions to $\X_k^*$. Hence,
$$ \limsup_k \mu(\{|x^*(g_k)|>\varepsilon\})>\varepsilon,$$
which contradicts the fact that $\{ g_k \}$ is scalarly
null in measure.

It  follows  from the claim  that there exists a 
subsequence $\{ g_{n_k} \}$ such that
 $$\underset \|x^*\|\le1 \to\sup\mu(\{|x^*(g_{n_k})|>2^{-k}\})<2^{-k}.$$
Clearly $\{ g_{n_k} \}$ is scalarly null a.e.   
\qed
\enddemo

However, we know of at least one space that fails ($D_\infty$).

\proclaim{Theorem 3.8} 
$C[0,1]$ fails property ($D_\infty$).
\endproclaim 

\demo{Proof}   
Consider a sequence $\{ f_n \}$ in  the  
unit ball of $L_\infty \left( C\left[0, 1\right]\right)$. 
With $\O_i = \left[0, 1\right]$, 
write 
$f_n(\cdot) \equiv f_n \left(\cdot, t\right) \: \O_1 \times \O_2 \to \Bbb R$ 
 so that 
\roster 
\item"{(i)}" $f_n(s, \cdot ) \: \O_2 \to \Bbb R$ \quad is in 
      $ C\left[0, 1\right]$ for almost all $s\in \O_1$ 
\item"{(ii)}" $\lav f_n (s,t) \rav \leqslant 1$ \quad  for each $t\in\O_2$ 
       for almost all $s\in \O_1$. 
\endroster
Such a sequence  
$\{ f_n \}$ is scalarly null in  measure  if and only if  
\roster
\item"{(iii)}"  $f_n ( \cdot, t ) \: \O_1 \to \Bbb R$ converges 
	 in measure to the null function for each $t\in\O_2$.  
\endroster
To see this, note that 
if $t\in\O_2$ is fixed, 
then  
$\left( x_t^* f_n \right) (\cdot) = f_n (\cdot, t )$ 
where  $x^*_t \in \left( C\left[0, 1\right]  \right)^*$ 
is the point evaluation at $t$ functional.  
As for the  reverse implication, 
assume that (iii) holds and let 
$ x^* \equiv \nu \in \left( C\left[0, 1\right]  \right)^*$  
be a finite regular positive Borel measure on 
$\Cal B\left(\O_2\right)$.  
It suffices to show that 
$  x^* f_n(\cdot) \equiv \int_{\O_2}  f_n (\cdot,t)  \, d\nu(t) $  
converges to the null function in $\mu$-measure.   
Towards this, let 
$\lambda = \mu \times \nu$ be the corresponding product measure 
on the completion $\Cal A$ of  $\Cal B\left(\O_1 \times \O_2\right)$. 
Then~(iii) implies that 
$f_n (\cdot, \cdot)\: \O_1 \times \O_2 \to\Bbb R$ 
converges to the null function  
in $\lambda$-measure and hence (by~(ii)) also in
$L_1\left( \O_1 \times \O_2, 
      \Cal A,  \lambda \right)$. 
Since 
$$
\int_{\O_1} \left( \int_{\O_2} \lav f_n (s,t) \rav \, d\nu(t) 
  \right) \, d\mu(s) ~=~ \iint_{ \O_1 \times \O_2} \lav f_n (s,t) \rav 
   \, d\lambda \rightarrow0  
$$ as $n\rightarrow\infty$, it follows that 
the sequence 
$\{ l_n \}$ of   $L_1$ functions  given by  
$$
  l_n(\cdot) ~\equiv~ \int_{\O_2} \lav f_n (\cdot,t) \rav \, d\nu(t)
$$ 
converges to the null function in $\mu$-measure, which gives the result.

For each positive integer $n$, 
let $\widetilde n$ be  its binary representation  as a finite 
sequence of $0$ and $1$'s. For     
$t\in \O$, let  
$t_3$ be its unique (nonterminating) ternary expansion 
into $0$, $1$, and $2$'s.   
For $1\leqslant k \leqslant n$, 
let $\Gamma(k,n)$ be the collection 
of all $k$-tuples  
$(n_1, n_2, \ldots, n_k)$  of positive integers  that satisfy   
$1\leqslant n_1 < n_2 <  \ldots <  n_{k-1} < n_k = n $. 
For $\gamma = (n_1, n_2, \ldots, n_k)$ in  $\Gamma(k,n)$,  
let $A_\gamma$ be the set of  
 $t\in\O$ for which  
$t_3$ is of the form 
$$ 0. ~\widetilde n_1 ~2~ \widetilde n_2 ~2~  \widetilde n_3 ~2~  
 \ldots \widetilde n_{k-1} ~2~  \widetilde n_k ~ 2 ~\ldots, 
$$ 
i\.e\. 
$$
A_{\gamma}  ~= ~(~ 
  0. ~\widetilde n_1 ~2~ \widetilde n_2 ~2~   
       \ldots \widetilde n_{k-1} ~2~  \widetilde n_k ~ 1~\overline{2} 
  \quad,\quad  
  0. ~\widetilde n_1 ~2~ \widetilde n_2 ~2~   
	\ldots \widetilde n_{k-1} ~2~  \widetilde n_k ~ \overline{2}  ~] \ .
$$   
For technical reasons, consider the subset
$$
\widetilde A_\gamma ~= ~(~ 
  0. ~\widetilde n_1 ~2~ \widetilde n_2 ~2~   
       \ldots \widetilde n_{k-1} ~2~  \widetilde n_k ~ 20~\overline{2} 
  \quad,\quad  
  0. ~\widetilde n_1 ~2~ \widetilde n_2 ~2~   
	\ldots \widetilde n_{k-1} ~2~  \widetilde n_k ~ 220\overline{2}  ~]  
$$
of $A_\gamma$ along with the corresponding unions 
$$ 
A^n_k ~=~ \bigcup_{\gamma\in\Gamma(k,n)} A_\gamma  
\text{\qquad and \qquad} 
\widetilde A^n_k ~=~ \bigcup_{\gamma\in\Gamma (k,n )} \widetilde A_\gamma  \ . 
$$       
The following properties of these sets will be used:             
\roster
\item $A^n_{k_1} \cap A^n_{k_2} ~=~ \emptyset$  \quad if    
      $k_1 \neq k_2$;  
\item  if $t \in \cap_j A^{n_j}_{k_j}$ for an increasing sequence $\{ n_j\}$, 
	then $\{ k_j \}$ is also (strictly) increasing;   
\item  if $\{ n_k \}_k$ is an increasing sequence of positive integers 
      and \newline $t_3 = 0. ~\widetilde n_1 ~2~ \widetilde n_2 ~2~ 
		     \widetilde n_3 ~2~  \ldots  $, 
      then $t \in \widetilde A^{n_k}_k$ for each $k$. 
\endroster 

For each admissible $n$ and $k$,    
find a continuous  function 
$g_k^n \: \Omega_2 \to [0,1]$  
that is supported on $A^n_k$ and takes the 
value $1$ on $\widetilde A_k^n$.   
Lexicographically order the dyadic interval 
$\{ I_k \}_{k\geqslant 1}$ of  $[0,1]$ 
and let $h_{k} \: \Omega_1 \to [0,1]$  
be the indicator function of  $I_k$.   
Define $f_n\: \O_1 \times \O_2 \to \Bbb R$   
by 
$$ 
     f_n (s,t) = \sum_{k=1}^{n} h_k(s)g^n_k (t) \ . 
$$
Clearly, the  corresponding sequence $\{ f_n \}$ is in the  
unit ball of $L_\infty \left( C\left[0, 1\right]\right)$ 
and it satisfies conditions (i),  (ii), and, by (1) and (2), 
also   (iii).  Thus $\{ f_n \}$ converges 
scalarly in measure to the null function.   

However, for any subsequence $\{ n_k \}_k$  of the 
positive integers, for  the 
corresponding  
point 
$t_3 = 0 . ~\widetilde n_1 ~2~ \widetilde n_2 ~2~ 
		    \widetilde n_3 ~2~  \ldots  $~, 
it follows from  
condition (3) that 
$f_{n_k}(s,t) = h_k (s)$, 
which does not go pointwise a\.e\. 
to the null function. 
\qed
\enddemo

We shall prove below (Corollary~4.11) that
$L_1$ fails ($D_p$) for  $1\leq p < \infty$,  
but we do not know what happens when $p=\infty$.

\proclaim{Question 3.9} Does $L_1$ enjoy ($D_\infty$)?  
\endproclaim 
This question can be reformulated  as a
question about functions of two variables as follows. 
For $n\ge1$, let $f_n(s,t)$ be real-valued functions on the unit
square which satisfy the following:
\roster
\item"(i)" $f_n(s,t)=\sum_{k=1}^{N(n)}1_{E_{n,k}}(s)g_{n,k}(t)$,
where $\{E_{n,k}\}_k$ is a partition of $[0,1]$ into    
sets of strictly positive measure;
\item"(ii)" $\int|g_{n,k}(t)|\,d\mu(t)\le1$ for all $1\le k \le N(n)$;
\item"(iii)" $F_n(s)=\int_A f_n(s,t)d\mu(t) \rightarrow 0$
in measure as $n\rightarrow\infty$ for every measurable set $A\subset
[0,1]$. \endroster
\proclaim{Question 3.9 paraphrased} For $\{f_n(s,t)\}$ as above,
does there  always exist a 
subsequence $\{f_{n_k}(s,t)\}$ such that
$\int_A f_{n_k}(s,t)d\mu(t) \rightarrow 0$ a.e. for each $A\in \Sigma$? 
\endproclaim
\heading{4.  $L_p(\X)$-BOUNDED  SEQUENCES}\endheading \vskip 10pt

We now investigate what happens when $L_\infty(\X)$-boundedness 
is weakened to $L_p(\X)$-boundedness.   
\subhead [4.i] From scalar convergence to weak a.e. convergence
\endsubhead
In this subsection, we shall use Dvoretzky's theorem on the existence
of almost spherical sections [Dv] to prove that for $1 \leq p <\infty$
none of the properties ($A_p$),  ($B_p$) nor
($C_p$) can hold in any infinite-dimensional Banach space.

Let us first recall the $q$-Pettis norm of an $L_0(\X)$ function 
(which might be infinite):    
$$ \|f\|_{\Cal P_q(\X)} = \sup_{x^*\in B(\X^*)}
\left( \int_\Omega |x^* f(\omega)|^q
\, d\mu\right)^{1/q} \  , $$
for $1 \leq q < \infty$. 
The building block used in our construction is
 the basic example of [DG] which we now recall
(and refine slightly) for the reader's convenience.
\proclaim{Proposition 4.1} Let $\X$ be an infinite-dimensional Banach space and
 let $E$ be a measurable subset of $\Omega$. Given 
$\varepsilon>0$, there exists $f\in L_\infty(\X)$ such that
$\|f(\cdot)\|=1_E(\cdot)$ and $\|f\|_{\Cal P_q(\X)}<2\varepsilon^{1/q}$
for each $1\le q<\infty$.
\endproclaim
\demo{Proof} 
Since $q \to \|f\|_{\Cal P_q(\X)}$ is an increasing function for a fixed 
$f \in L_0 (\X)$, it suffices to consider only 
$2\le q<\infty$.  
First we prove the result for $E=\Omega$.

Let
$$
\{I^n_k=[(k-1)/2^n,k/2^n): n\ge0 
\text{~and~} 1\le k \le 2^n\}
$$
be the collection of
 dyadic subintervals of $\Omega$. By Dvoretzky's Theorem there exist
unit vectors $\{e^n_k\}_{k=1}^{2^n}
$ in $\X$ such that
$$ \frac{1}{2} \left( \sum_{k=1}^{2^n} |a^n_k|^2\right)^{1/2}
\le \left\|\sum_{k=1}^{2^n} a^n_ke^n_k \right\| 
\le 2\left(\sum_{k=1}^{2^n} |a^n_k|^2\right)^{1/2}\tag1$$
for all real numbers $a^n_k$.
Define $f_n:\Omega\rightarrow\X$ by
$$f_n(\omega)=\sum_{k=1}^{2^n}1_{I^n_k}(\omega)e^n_k.$$
Note that $\|f_n(\omega)\| =1$ for all $\omega\in \Omega$.
 Fix $x^*\in B(\X^*)$. Then (1) implies that $\left(\sum_{k=1}^{2^n}
|x^*(e^n_k)|^2\right)^{1/2} \le2.$ Thus, for $q\ge2$, we have
$$
\align \int_\Omega \left|x^*f_n(\omega)\right|^q\,d\mu 
&= 
\int_\Omega \left|\sum_{k=1}^{2^n} 
  x^*(e^n_k)1_{I^n_k}(\omega)\right|^q\,d\mu \\
&= 
2^{-n}\sum_{k=1}^{2^n}|x^*(e^n_k)|^q\\
&\le 
2^{-n}\left(\sum_{k=1}^{2^n}|x^*(e^n_k)|^2\right)^{q/2} \le 2^{q-n}~,
\endalign
$$
and so 
$\|f_n\|_{\Cal P_q(\X)}\le 2 \cdot \left(2^{-n}\right)^{\frac{1}{q}}$, 
which gives the result. 

An analogous
construction  can be carried out in any set $E$ of positive measure,
and the result is  trivial anyhow for a set $E$ of measure zero. 
\qed\enddemo
Now we beef up the previous result.
\proclaim{Proposition 4.2} Let $\X$ be an infinite-dimensional Banach space and
 let $h$ be a non-negative countably-valued measurable
 function defined on $\Omega$. Given $\varepsilon>0$ 
and $1 \leq q_0<\infty$
 there 
exists $f\in L_0(\X)$ with the following properties:
\roster
\item
  $\|f(\cdot)\|=h(\cdot)$;
\item
$\|f\|_{\Cal P_q(\X)} < \infty$ for each $1 \leq q<\infty$;
\item  
$\|f\|_{\Cal P_{q_0}(\X)}<\varepsilon$. 
\endroster 
\endproclaim
\demo{Proof}
 Write $h=\sum_{k=1}^\infty a_k 1_{E_k}$, where the $a_k$'s are
 positive numbers and the $E_k$'s are disjoint measurable sets.
Select positive numbers $\{ \varepsilon_k\}$ 
such that 
$\sum_{k=1}^\infty a_k^q \varepsilon_k$ 
is finite for each $ 1 \leq  q<\infty$ and 
$\sum_{k=1}^\infty a_k^{q_0}
\varepsilon_k < \left({\varepsilon}/{2}\right)^{q_0}$. 
 By Proposition~4\.1, for each $k$ there exists
$f_k \in L_\infty(\X)$ such that $\|f_k(\cdot)\|= 1_{E_k}(\cdot)$ 
and $\|f_k\|_{\Cal P_q(\X)} < 2\varepsilon_k^{1/q}$. Clearly,
$f= \sum_{k=1}^\infty a_k f_k$ has the required properties.\qed
\enddemo  

 \proclaim{Theorem 4.3} 
Let $\X$ be an infinite-dimensional Banach space
and
 let $g$ be any non-negative measurable function that 
 is {\it not} essentially-bounded. There exists a sequence
$\{g_n\}$ in $L_0(\X)$ such that the following hold:
\roster
\item $$\mu(\{\omega:\|g_n(\omega)\|>t\})\le \mu(\{\omega:g(\omega)>t\})$$
for all $n\ge1$ and for all $t>0$;
\item $$\sum_{n=1}^\infty \|g_n\|_{\Cal P_q(\X)}<\infty$$
for each $ 1 \leq q<\infty$;
\item  $\{g_n\}$ converges scalarly a\.e\. to the null function;
\item for each subsequence $\{g_{n_j}\}$ there exists a set
$A\subset \Omega$ of full measure such that
$$\limsup_j \|g_{n_j}(\omega)\|=\infty $$ 
for each  $\omega\in A$. 
\endroster
In particular, no subsequence of $\{g_n\}$ converges weakly on
any set of strictly positive measure. \endproclaim
\demo{Proof}  
Let $h$ be a non-negative countably-valued  measurable function on
$\Omega$ which is {\it not} essentially bounded and which satisfies
$h(\omega)\le g(\omega)$ for $\omega\in \Omega$. Use Proposition~4\.2 to
construct
 a sequence $\{g_n\}$ of {\it independent} $\X$-valued
random variables such that 
\roster
\item"{(i)}" each $\| g_n \|$ has the
same distribution as $h$, 
\item"{(ii)}"  $\|g_n\|_{\Cal P_q(\X)}$ is finite for each 
$ 1 \leq q<\infty$, and
\item"{(iii)}"  $\|g_n\|_{\Cal P_{n}(\X)}\le 2^{-n}$ \ .   
\endroster 
 Clearly (1) is satisfied.  
Condition~(2) follows from  the observation that,  
if $N \in \Bbb N$ and $1 \leq q \leq N$, then by (ii) and (iii)   
$$\align 
\sum_{n=1}^\infty \| g_n\|_{\Cal P_q(\X)} 
\sum_{n=1}^N
\|g_n\|_{\Cal P_q(\X)} +   \sum_{n=N+1}^\infty \|g_n\|_{\Cal P_n(\X)}\\
 \sum_{n=1}^N \|g_n\|_{\Cal P_q(\X)} +
 \sum_{n=N+1}^\infty 2^{-n} 
\quad < \quad \infty~. \endalign
$$ 
Clearly, (3) follows from (2) using Fact~1\.5. To prove (4),
fix a subsequence $\{g_{n_j}\}$. Then, for each $M>0$,    
$$\sum_{j=1}^\infty \mu(\{\omega:\|g_{n_j}(\omega)\|>M\})
=\sum_{j=1}^\infty \mu(\{\omega: h (\omega)>M\})=\infty$$
since  $h$  does not
belong to $L_\infty$. So by the Borel-Cantelli lemma
$\|g_{n_j}(\omega)\|>M$ infinitely often a\.e\. \qed
\enddemo
An appropriate choice of the measurable function $g$
(e\.g\. $g(\omega)= |\log \omega|$) in Theorem~4.3
yields the following corollary.
\proclaim{Corollary 4.4} Let $\X$ be a Banach space
and let $1\le p<\infty$.  The following
are equivalent:
\roster
\item $\X$ is finite-dimensional;
\item $\X$ satisfies $(A^\prime_p)$;
\item $\X$ satisfies $(A_p)$;
\item $\X$ satisfies $(B_p)$;
\item $\X$ satisfies $(C_p)$.
\endroster \endproclaim
\demo{Remarks} 1. Theorem~4.3 shows that there is no analogue
for scalar convergence of the {\it uniform boundedness principle}:
if $\X$ is infinite-dimensional then `scalar boundedness a\.e\.'
does {\it not} imply `norm-boundedness a\.e\.' \newline
2. If $\|f_n\|_{L_q(\X)}\rightarrow0$ then clearly some subsequence
is $\X$-norm-null a\.e. However, condition (2) of Theorem 4\.3 suggests that
 searching for a non-trivial
scalar integrability condition
which implies weak a\.e\. convergence is probably futile. \enddemo

\subhead [4.ii] From scalar convergence in measure or in $L_p$
to scalar a.e\. convergence    \endsubhead      
In this subsection we examine 
the properties ($D_p$) and ($E_p$) more closely for $1 \leq p < \infty$.   

First 
we recall some notation from [Pr].
Let $\{ x_n \}$ be a basic sequence in a Banach space $\X$
with coefficient functional sequence $\{ x_n^* \}$ in $\X^*$.
 A family
$\{ \X_n \}$ of finite-dimensional subspaces of $[x_n]$ is a
{\it blocking} of $\{ x_n \}$ provided there exists an increasing
sequence of integers $\{ n_k \}$ with $n_1=1$ such that
$\X_k=[x_i]_{i=n_k}^{n_{k+1}-1}$ for each~$k$.   
For $1\le p \le\infty$, 
if there is a positive constant $c$ so that, 
 for each collection of vectors $\{y_i\}_{i=1}^n$, where $y_i \in \X_i$, 
$$ \ln \sum_{i=1}^n y_i \rn ~\leq ~
   c \ln ~ \left( \ln y_i \rn \right) ~ \rn_{\ell_p} \ ,  
$$
resp. 
$$ c \ln ~ \left( \ln y_i \rn \right) ~ \rn_{\ell_q} ~\leq ~
  \ln \sum_{i=1}^n y_i \rn  \ , 
$$
then we say that the blocking $\{ \X_k\}$ satisfies an 
upper  (resp. lower) $p$-estimate.

\proclaim{Theorem 4.5} 
Fix $1 < q \le \infty$ and let $q^{\prime}$ be its 
conjugate exponent.  
Suppose that $\X$ has a basis $\{ x_n \}$ with the 
property that each blocking of this basis     
satisfies an upper $q$-estimate.  Then $\X$ 
enjoys ($D_{p}$) for each $q^{\prime} \leq p \leq \infty$. 
\endproclaim 

\demo{Proof} 
Wlog, $p = q^{\prime}$.    
Fix a sequence $\{ f_n \}$  
in $B(L_{p}(\X))$
that is scalarly null in measure.  We need to extract 
an scalarly null a.e.~subsequence.  To this end, 
let $\{ \X_k \}$, $\{ P_k \}$, $\{ g_k \}$, and $\{ h_k \}$ be as provided 
from Lemma~3.7.  It suffices to show that 
$\{ g_k \}$ converges to the null function 
scalarly a\.e.  

Fix $x^* \in \X^*$ and let $x^*_k = x^* \circ P_k \in X^*$. 
Note that 
$$ 
\lav x^* g_k (\omega) \rav \leq \ln x^*_k \rn ~ \ln g_k(\omega) \rn_{\X} \ . 
$$ 
Thus, for $\e > 0$ fixed 
$$ 
  \mu (\{ \lav x^* g_k \rav \geq \e \} )  
  ~\leq~ 
   \left[ \frac{\ln x^*_k \rn}{\e} ~ \ln g_k \rn_{L_{p}} 
   \right]^{p} \ . 
$$ 
Note that each $\ln g_k \rn_{L_{p}}$ is bounded above by $2K$ 
where $K$ is the basis constant of $\{ x_n \}$. 
Thus 
$$
\sum_{k=1}^\infty  \mu (\{ \lav x^* g_k \rav \geq \e \} )  
$$ 

Since the blocking $\{ \X_k \}$ satisfies an upper $q$-estimate  
(say with constant $C$),  
$$ \ln ~ \left( \ln x^*_k \rn \right) ~ \rn_{\ell_p} ~\leq ~ 2CK
 \liminf_n  \ln \sum_{k=1}^n x^*_k \rn   ~\leq ~2 CK^2 \ln x^* \rn \ . 
$$ 
Thus 
$ 
\sum_k  \mu (\{ \lav x^* g_k \rav \geq \e \} ) < \infty$. 
So by Borel-Cantelli,  
$\{ x^* g_k \}$ converges to the null function, as needed. 
\qed 
\enddemo

Minor variations in the above proof 
 give that  
Theorem~4.5 remains true   if the word {\it basis} is replaced 
by {\it finite-dimensional decomposition\/}.  

However, there are many spaces that  fail ($E_p$) 
(and hence fail ($D_p$)).
\proclaim{Theorem 4.6} 
Fix $1 < q \leq \infty$ and let $q^\prime$ be 
its 
conjugate exponent. 
Suppose that $\X$ contains a  weakly null semi-normalized basic
sequence  $\{ x_n \}$ which  satisfies a lower $q$-estimate. 
Then $\X$ fails ($E_p$) for each $1\le p < q\prime$.  
\endproclaim
\demo{Proof} We may assume without loss 
of generality that $\X=[x_n]$.
Fix $p\in [1,q^\prime)$ and choose
$q_0 \in \left(\frac{1}{q'} , \frac{1}{p} \right)$. 
Let $\{ g_n \}$ be a sequence of i.i.d. 
random variables defined on $\Omega$ 
with the same distribution as 
$  g_0 (t) = t^{-q_0}   $. 
Define  
$f_n \: \Omega \to \X$ by 
$$ 
   f_n (\cdot)~ =  ~g_n(\cdot)~ x_n  \ . 
$$ 

Since $g_0 \in L_p$ and since $ \{ x_n \}$ is semi-normalized, 
$\{ f_n \}$ is an $\lpx$-bounded sequence.   
For each $x^* \in \X^*$, 
$$ 
  \ln x^* f_n \rn_{L_p} ~=~ \lav x^* (x_n) \rav ~  
   \ln g_n \rn_{L_p} \ . 
$$ 
Thus $\{f_n \}$ converges scalarly in $L_p$ to the null function.  

Fix a subsequence $\{ f_{n_j} \}$ of $\{f_n \}$. 
It suffices to show that $\{ f_{n_j} \}$  is not 
scalarly null a\.e\. To this end,  
let $\{ x^*_n \}$ be the sequence of biorthogonal 
functionals
satisfying $x^*_n (x_m) = \delta_{nm}$.  
Since $\{ x_n \}$ is semi-normalized and 
 satisfies a lower $q$-estimate, it follows that 
$\{ x^*_n \}$ satisfies   an upper $q^\prime$-estimate.   
Consider the element $x^*\in \X^*$ given by 
$x^* = \sum_j j^{-q_0}~x^*_{n_j}$, which  
 converges in $\X^*$
since $\{ x_n^*\}$ satisfies an upper $q^\prime$-estimate.   
Fix $M>0$. Since  
$$
\mu\{ x^* f_{n_j} > M \} ~=~ 
\mu\{ j^{-q_0} g_{n_j} > M\} 
M^{-\frac{1}{q_0}} ~j^{-1} \ , 
$$  
we see that $\sum_j \mu\{ x^* f_{n_j} > M \} = \infty$, and 
so by the Borel-Cantelli lemma 
there is a set $A$ of full measure such that 
if $\omega \in A$ then 
$|x^* f_{n_j} (\omega)| > M$ infinitely often.  
Thus this subsequence does not 
converge scalarly a\.e\.   
\qed
\enddemo
By a theorem of Prus \cite{Pr}
every {\it nearly uniformly convex} space
 (see \cite{H} for the definition
of this property)
satisfies the
hypothesis of Theorem~4.6 for some $1 < q < \infty$  and so we obtain the
following corollary. 
\proclaim{Corollary 4.7} Suppose that $\X$ contains a super-reflexive
or (more generally)  a nearly uniformly convexifiable 
infinite-dimensional  subspace. 
Then $\X$ fails ($E_p$) for 
some $1 <  p < \infty$. \endproclaim

The following two corollaries follow  from the previous 
two theorems  by  considering the standard unit vector basis 
of $\ell_p$ and $c_0$.

\proclaim{Corollary 4.8} 
Fix $1 < q < \infty$ and  let $q^\prime$ be its conjugate exponent. 
Then
$\ell_q$ has ($D_p$) (resp.~($E_p$)) 
if and only if  \, $1 <  q^\prime \leq p \leq \infty$. 
\endproclaim

\proclaim{Corollary 4.9}
$c_0$ has ($D_p$) (and thus ($E_p$)) 
for each $1 \leq p \leq \infty$.
\endproclaim 
Note that in Corollary~4.8,  
if  $q\downarrow1$ 
then $p\uparrow\infty$,  
which suggests that $\ell_1$ should fail $D_p$ for  
each $1 \leq p<\infty$.
Strangely, however,  the truth  is the complete opposite as
was proved in Theorem~3.6 above: $\ell_1$ has $D_p$ for all $1 \le p
\le \infty$.

Finally, we determine the range of values of $p$ for which
$L_q$ satisfies $(D_p)$. 
\proclaim{Corollary 4.10} 
Fix $1 < q < \infty$ and  let $q^\prime$ be its conjugate exponent. Then 
$L_q$ has ($D_p$) (resp.~($E_p$))    
if and only if \, $ \max{(2, q^\prime)} \leq p \leq \infty$. 
\endproclaim

\demo{Proof} 
Since $\ell_2 $  and $\ell_q$ each embed  into $L_q$, 
it follows from Corollary~4.8 that $L_q$ fails  ($E_p$) 
for $1\leq p < \max{(2, q^\prime)}$. 
Since $L_q$ has type $\min({2, q})$  
and the Haar system $\Cal H$ forms an 
unconditional basis for $L_q$,  each blocking of $\Cal H$ satisfies 
an upper $\min({2, q})$-estimate.   So by Theorem~4.5,  
if   $p = \max{(2, q^\prime)}$ then 
$L_q$ enjoys ($D_{p})$.  
\qed
\enddemo
Let $H_1$ denote the Hardy space of analytic functions on the
unit disk in the complex plane with the usual $L_1$ norm 
(see e\.g\. \cite{Ru}). It is known that $H_1$ contains
subspaces that are isomorphic to $\ell_q$ for $1\le q \le 2$
(see e\.g\. \cite{Di}).
Hence
Corollary~4.8 implies the next result (cf. Question~3.9).

\proclaim{Corollary 4.11} 
$H_1$ (and therefore also $L_1$) fails ($E_p$) for $1 \le p <  \infty$. 
\endproclaim

\demo{Remark} We do not know of a reflexive 
space that satisfies ($D_1$).
However, we resist making the obvious conjecture.
\enddemo

\head 5. Completeness \endhead
In this section we prove some completeness results 
for the topologies of scalar convergence considered
in this paper. First we recall the appropriate definitions. 
Let $E$ be a topological vector space. A sequence $\{x_n\}$ in $E$
is a Cauchy sequence  if for every zero-neighborhood $U$
there exists $N\ge1$ such that $x_n-x_m\in U$ for all
$n,m \ge N$. We shall say that $E$ is complete if every
Cauchy sequence converges.

\proclaim{Theorem 5.1} Let $\X$ be an infinite-dimensional Banach space.
The topologies of scalar convergence in measure and scalar
convergence in $L_p$ (\/$1\le p <\infty$) are incomplete.
\endproclaim
\demo{Proof} This result is very similar in spirit to the fact that
the Pettis norm is incomplete whenever $\X$ is infinite-dimensional
\cite{JK}. We refer the reader to the proof of the incompleteness
of the Pettis norm that is given 
in \cite{DG}. The construction there, 
which utilizes Dvoretzky's Theorem on almost spherical sections \cite{Dv},
can easily be 
modified, using the estimates of Proposition~4.1, 
to construct a sequence of functions that is 
Cauchy but not convergent in the topology of scalarly 
convergence in measure  and  $L_p$ for $1 \leq p < \infty$.  \qed
\enddemo
\proclaim{Theorem 5.2} The topology of scalar convergence 
in $L_\infty$
is complete if and only if $\X$ is weakly sequentially complete.
\endproclaim
\demo{Proof} First suppose that $\X$ is not weakly sequentially  
complete. Let $\{x_n\}$ be a weak Cauchy sequence that does not
converge weakly 
and let $f_n(\omega)=x_n$ ($n\ge1$). Clearly, $\{f_n\}$
is a Cauchy sequence in the topology of scalar convergence in
$L_\infty$.  
By Proposition~3.3,    a   
limit of this sequence, say $f$,  would have to satisfy
$f(\omega)=\text{ weak-lim } f_n(\omega)$ almost everywhere. 
 Hence
$\{f_n\}$ does not converge. 

For the converse, suppose that $\X$
is weakly sequentially complete. Let $\{f_n\}$ be a Cauchy
sequence in the topology of scalar convergence
in $L_\infty$. By adapting the proof of Proposition~3.3,   
we see that the weak sequential completeness of $\X$ 
guarantees that  there exists a function $f:\Omega\rightarrow \X$
such that  $f(\omega)=\text{weak-lim }f_n(\omega)$ a.e.
By the Pettis measurability theorem \cite{P} and Proposition~1.1,
$f$ is strongly-measurable, i\.e\., $f\in L_0(\X)$. It now 
follows easily from the fact that $\{f_n\}$ is a Cauchy sequence 
that $\{f_n\}$ converges to $f$ scalarly in $L_\infty$. \qed 
\enddemo
Of more relevance to this paper  is the convergence
of a pointwise-bounded or an $L_p(\X)$-bounded Cauchy sequence.
We investigate this question next for the topology of
scalar convergence in measure. For brevity's sake we
shall say that a sequence is ``scalarly Cauchy in measure"
if it is a Cauchy sequence for the topology of scalar
convergence in measure.
\proclaim{Theorem 5.3} Let $\X$ be a reflexive
Banach space. Then each pointwise-bounded sequence $\{f_n\}$ in
$L_0(\X)$ that is scalarly Cauchy in measure   
converges  scalarly  in measure.
\endproclaim 
\demo{Proof} We may assume that $\X$ is separable and hence that
$\X^*$ is separable. Arguing now as in Proposition~3.1 there 
exists a subsequence $\{f_{n_k}\}$ such that $\{f_{n_k}(\omega)\}$
is a weakly Cauchy sequence in $\X$ almost everywhere.  
 Since $\X$ is reflexive
it is weakly sequentially complete and so (by the Pettis Measurability Theorem 
and Proposition~1.1) 
there exists $f$ in $  L_0(\X)$ 
 such that $f_{n_k}$ converges to $f$ 
weakly a\.e\., thus also 
scalarly in measure, which is enough.  \qed
\enddemo

\proclaim{Theorem 5.4} 
Let $\X$ be a reflexive
Banach space.  Then each  $L_1(\X)$-bounded sequence $\{f_n\}$ 
that    is  scalarly Cauchy  in measure   
converges  scalarly  in measure.  
\endproclaim 
\demo{Proof} First, we may assume that $\X$ is separable. By a deep 
result of Zippin \cite{Z} every separable reflexive 
Banach space is isomorphic
to a closed subspace of a reflexive Banach space with a basis.
So we may assume that $\X$ is isomorphically embedded into a reflexive Banach 
space $\Y$ with a basis. Clearly, $\{f_n\}$
is scalarly Cauchy in measure when viewed as a sequence in $L_0(\Y)$. 
By Proposition~1.1, it suffices to show that 
 $\{f_n\}$ converges to some $f$ in $  L_0(\Y)$.

Since $\Y$ is reflexive, a normalized basis $\{e_k\}$ for $\Y$  
is both
boundedly complete and shrinking \cite{LiT}. Let $C$ be the 
basis constant of $\{e_k\}$.
For each $n$, we can
expand $f_n$ with respect to the basis $\{e_k\}$ thus:
$$f_n(\omega)=\sum_k f_{n,k}(\omega)e_k.$$
For each $k$, the sequence $\{f_{n,k}\}_n$ is  Cauchy  in 
measure, and hence converges in
measure to some $g_k\in L_0$. 
By Fatou's Lemma, we have
$$\align \int_\Omega \sup_N \left\|\sum_{k=1}^N g_k(\omega)e_k\right\|
\,d\mu
&\le C \int_\Omega \liminf_{N\rightarrow\infty}  
\left\|\sum_{k=1}^N g_k(\omega)e_k\right\|\,d\mu\\
&\le C \liminf_{N\rightarrow\infty} \int_\Omega 
\left\|\sum_{k=1}^N g_k(\omega)e_k\right\|\,d\mu\\
&\le C \liminf_{N\rightarrow\infty} 
\left(\liminf_{n\rightarrow\infty}\int_\Omega 
\left\|\sum_{k=1}^N f_{n,k}
(\omega)e_k\right\|\,d\mu\right)\\
&\le C \liminf_{N\rightarrow\infty} \left(C\liminf_{n\rightarrow\infty} 
\int_\Omega \left\|\sum_{k=1}^\infty 
f_{n,k}(\omega)e_k\right\|\,d\mu\right)\\
&\le C^2 \sup_n \|f_n\|_{L_1(\Y)}<\infty. \tag1 \endalign$$
Hence
$$\sup_N\left\|\sum_{k=1}^N g_k(\omega)e_k\right\|<\infty\quad a.e.$$
Since $\{e_n\}$ is boundedly complete it follows that
$f(\cdot) \equiv \sum_{k=1}^\infty g_k(\cdot)e_k$ is 
in $L_0( \X)$;  
 moreover, it follows from (1) that $f \in L_1(\Y)$.
 
Fix $y^*\in \Y^*$ and $N\ge1$. Clearly,
$$y^*\left(\sum_{k=1}^N f_{n,k}e_k\right)\quad \rightarrow \quad
y^*\left(\sum_{k=1}^N g_ke_k\right)\tag2$$
in measure as $n\rightarrow\infty$.
Let $\alpha_n$ denote the norm of the restriction of $y^*$
to $[e_k]_{k\ge n}$. Since $\{e_n\}$ is a shrinking basis, 
 $\alpha_n\rightarrow0$ as $n\rightarrow\infty$. Now
$$\align
\int \left|y^*\left(\sum_{k=N+1}^\infty f_{n,k}e_k\right)\right|\,d\mu
&\le \alpha_{N+1} \int \left\|\sum_{k=N+1}^\infty f_{n,k}e_k\right\|\,d\mu\\
&\le \alpha_{N+1}(1+C)\sup_n \|f_n\|_{L_1(\Y)}\quad \rightarrow \quad 0\tag3
\endalign$$
as $N \rightarrow \infty$.
Combining (2) and (3) we see that $\{f_n\}$ 
converges to $f$ scalarly in measure.   
\qed
\enddemo

The proof of Theorem~5.3 apparently uses only
the weak sequential completeness of $\X$ and the fact that 
$\X^*$ has the RNP.  
However, by  Rosenthal's $\ell_1$ theorem [Ro],  
   these two
properties are {\it equivalent} to  $\X$ being reflexive.    
Clearly, a necessary condition
for the conclusion of Theorems~5.3 and~5.4 
to hold is 
 that $\X$ is weakly sequentially complete,    
and 
when $\X$ has an unconditional basis this condition is also sufficient, 
as our next two results show.  
However, we have not been able to determine 
general necessary and   sufficient conditions 
on $\X$ so that the 
conclusions of Theorems~5.3 and~5.4 
 hold. 
In view of the next two theorems, which establish 
the desired conclusions for $\ell_1$, it is clear that the
method of proof of 
Theorem~3.2 will not be of use in this situation.

\proclaim{Theorem 5.5}  
Let $\X$ be a weakly sequentially complete
Banach space with an unconditional basis.  
Then each pointwise-bounded sequence $\{ f_n \}$ in 
$L_0(\X)$ that is scalarly Cauchy in
measure   converges  scalarly in measure.   
\endproclaim 
\demo{Proof} Let $\{e_n\}$ be a normalized unconditional basis for $\X$. 
We may assume, without loss of generality, that
$$\left\|\sum_n a_n e_n\right\| \le \left\|\sum_n b_n e_n\right\|\tag1$$ 
whenever
$|a_n| \le |b_n|$ for all $n$. The fact that $\X$ is 
weakly sequentially complete implies that $\{e_n\}$ is
boundedly complete \cite{LiT}. 

By assumption, 
$$\sup_n \|f_n(\omega)\|=M(\omega)<\infty \quad \text{a.e.} \ . \tag2$$
Also, for each $n$, we can
expand $f_n$ with respect to the basis $\{e_k\}$ thus:
$$f_n(\omega)=\sum_k f_{n,k}(\omega)e_k \ .$$
For each $k$, the sequence $\{f_{n,k}\}_n$ is a Cauchy sequence in
the topology of convergence in measure, and hence converges in
measure to some $g_k$. Now (1) and (2) imply that 
$$\sup_N \left\|\sum_{k=1}^N g_k(\omega)\right\|\le M(\omega) \quad  \text{a.e.} \ . $$  
Since $\{e_n\}$ is boundedly complete it follows that
$f(\cdot)=\sum_{n=1}^\infty g_n(\cdot)e_n$ is in $L_0(\X)$. 

Let $h_n = f-f_n$, and so 
$$
 h_n(\cdot)= \sum_k ( g_k-f_{n,k} ) (\cdot)~e_k  \ . 
$$        
 To complete the proof of the theorem,  it suffices to
show that $\{h_n\}$  
is scalarly null in measure. 

So suppose,
to derive a contradiction, that  $\{h_n\}$ is not scalarly null in measure.
Then there exists $x^*\in S(\X^*)$ and $\varepsilon>0$ such that 
$$ \mu \{ \omega: |x^*h_n(\omega)| > \varepsilon\}> \varepsilon\tag3$$
for infinitely many $n$. 
  
The gliding hump argument of Lemma~3.7   yields a subsequence
$\{h_{n_k}\}_k$  and a blocking $\{\X_k\}$ of the basis such that
each $h_{n_k}$ satisfies~(3) and 
$$\mu\{\omega: \|h_{n_k}(\omega) - P_kh_{n_k}(\omega)\|>
\varepsilon/4 \} <\varepsilon/4,\tag4$$
where $P_k$ is the natural projection of $\X$ onto $\X_k$. We may
define $y^*\in \X^*$ by defining its action on each $x_k\in\X_k$:
$$
  y^*(x_k)=0\quad\text{($k$ odd)};\quad y^*(x_k)=x^*(x_k)
  \quad\text{($k$ even)} \ . \tag5 
$$ 
Then by (1), we have $\|y^*\|\le\|x^*\|\le1$, and so from (3),
 (4) and (5) we deduce that
 $$ \mu \{ \omega: |y^*h_{n_k}(\omega)| > \varepsilon/4\}< \varepsilon/4
 \quad\text{($k$ odd)}\tag6$$
 while
 $$ \mu \{ \omega: |y^*h_{n_k}(\omega)| > \varepsilon/2\}> 3\varepsilon/4
 \quad\text{($k$ even)}.\tag7$$
 Clearly, (6) and (7) contradict the fact that $\{h_n\}$ is scalarly
 Cauchy  in measure \qed
 \enddemo
 \proclaim{Theorem 5.6} 
Let $\X$ be a weakly sequentially complete
Banach space with an unconditional basis.  
Then each $L_1(\X)$-bounded sequence $ \{ f_n \}$ 
that is  scalarly Cauchy in
measure  converges scalarly  in measure.
\endproclaim 
\demo{Proof} Let $\{e_k\}$ be a normalized unconditional
basis for $\X$.  Let
$$f_n(\omega)=\sum_k f_{n,k}(\omega)e_k$$
 be the expansion of $f_n$. 
Now, arguing as in the first half of Theorem~5.4,
it can be shown that, for each $k$, $f_{n,k}\rightarrow g_k$ 
in measure as
$n \rightarrow \infty$ and that $f(\cdot)\equiv\sum_k g_k(\cdot) e_k$
belongs to $L_1(\X)$. Now, arguing as in second half of 
Theorem~5.5, one uses the
unconditionality of $\{e_k\}$ to prove that $\{f_n\}$ converges to
$f$ scalarly in measure. \qed \enddemo

\demo{Remark} Note that Theorems 5.5 and 5.6 apply to both $\ell_1$
and $H_1$. \enddemo

Finally, straightforward modifications to the proofs of Theorems~5.4
and 5.6 yield the following.
\proclaim{Theorem 5.7} Fix $1\le p< \infty$. Let
$\X$ be a weakly sequentially complete 
Banach space that is either reflexive or has an unconditional basis.
Then each $L_p(\X)$-bounded sequence that is scalarly Cauchy
in $L_p$  converges scalarly in $L_p$.
\endproclaim

\widestnumber\no{[GGG]Z}
\def\n #1{\no{[\bf #1]}}


\Refs

\ref\n{C}
\by       Kai Lai Chung 
\book     A Course in Probability  Theory 
\bookinfo Probability and Mathematical Statistics 
\vol    21
\publ     Academic Press 
\publaddr New York 
\yr 1974
\endref
\ref\n{DJ}   
\by        W. J. Davis and W. B. Johnson   
\paper     Weakly Convergent Sequences of Banach Space Valued 
	   Random Variables 
\inbook    Banach Spaces of Analytic Functions
\eds       J. Baker, C. Cleaver, and J. Diestel
\bookinfo  Lecture Notes in Math.
\vol       604
\publ      Springer-Verlag 
\publaddr  New York-Berlin 
\pages     29--31 
\yr        1977
\endref

\ref\n{DU}
\by         J. Diestel and J.J. Uhl, Jr.
\book       Vector Measures
\bookinfo   Math. Surveys, no. 15
\publ       Amer.\ Math.\ Soc.
\publaddr   Providence, R.I. 
\yr         1977
\endref


\ref \n{Di} \paper Intersection of Lebesgue spaces $L_1$ and $L_2$
\jour Proc. Amer. Math. Soc. \vol 103 \yr 1988 \pages 1185--1188 
\by S.J. Dilworth \endref

\ref\n{DG}
\by       S.J. Dilworth and Maria Girardi 
\paper    Bochner vs. Pettis norms\rom: examples and results
\inbook   Banach Spaces   
\bookinfo Contemp. Math. 
\ed       Bor-Luh Lin and William B. Johnson
\vol      144
\publ     American Mathematical Society 
\publaddr Providence, Rhode Island 
\yr       1993 
\pages    69--80
\endref


\ref\n{Dv}
\by      A. Dvoretzky
\paper   Some results on convex bodies and Banach spaces
\inbook  Proc. Internat. Symp. on Linear Spaces
\publaddr  Jerusalem 
\yr      1961
\pages   123--160
\endref
 






\ref\n{E} 
\by      G. A. Edgar 
\paper   Asplund Operators and A. E. Convergence  
\jour    J. Multivariate Anal. 
\vol     10
\yr      1980 
\pages   460--466 
\endref 

\ref\n{GS}
\by      N. Ghoussoub and P. Saab
\paper   Weak Compactness in Spaces of 
	 Bochner Integrable Functions and the Radon-Nikod\'ym property
\jour    Pacific J.  Math.
\vol     110
\issue   1
\yr      1984
\pages   65--70
\endref

\ref \n{H} \by R. Huff \paper Banach spaces which are nearly
 uniformly convex
\jour Rocky Mountain J. Math. \vol 10 \yr 1980 \pages 743--749 \endref

\ref\n{J} \by R. C. James \paper Super-reflexive spaces with bases
\jour Pacific J. Math \vol 41 \yr 1972 \pages 409--419 \endref
\ref\n{JK}
\by     Liliana Janicka and Nigel J. Kalton
\paper  Vector Measures of Infinite Variation
\jour   Bull. Polish Acad. Sci. Math.
\vol    XXV
\issue    3
\yr    1977
\pages 239--241
\endref

\ref\n{WBJ} 
\by W.B.~Johnson
\finalinfo Personal Communication   
\endref 

\ref\n{LeT} 
\by       Michel Ledoux and Michel Talagrand  
\book     Probability in Banach Spaces, Isoperimetry and Processes  
\bookinfo Ergebnisse der Mathematik und ihrer Grenzgebiete  
\publ     Springer-Verlag 
\publaddr Berlin 
\yr       1991
\endref 
\ref\n{LiT}
\by J. Lindenstrauss and L. Tzafriri
\book Classical Banach Spaces I: Sequence Spaces 
\publ Springer-Verlag   \publaddr Berlin-Heidelberg
\yr 1977
\endref


\ref \n{P}  
\by     B.J. Pettis  
\paper  On integration in vector spaces  
\jour   Trans. Amer. Math. Soc.  
\vol    44
\yr     1938 
\pages  277--304  
\endref  

\ref\n{Pr} \by Stanislaw Prus \paper Nearly uniformly smooth Banach spaces
\jour Boll. Un. Mat. Ital. \vol 7 \yr 1989 \pages 507--521 \endref
\ref\n{Ro} \by H.P. Rosenthal \paper A characterization of
Banach spaces containing $\ell_1$ \jour Proc. Nat. Acad.
Sci. U.S.A. \vol 71 \pages 2411--2413  \yr 1974
\endref

\ref\n{Ru} 
\by        Walter Rudin 
\book      Real and complex Analysis  
\publ      McGraw-Hill 
\publaddr  New York 
\yr        1974
\endref


\ref\n{S} 
\by Charles Stegall 
\paper The Radon-Nikod\'ym Property in Conjugate Banach Spaces 
\jour   Trans. Amer. Math. Soc.  
\vol    206 
\yr   1975
\pages  213--223
\endref

\ref\n{U}
\by J. J. Uhl, Jr. 
\paper A note on the Radon-Nikod\'ym property for Banach spaces 
\jour Rev. Roumaine Math. Pures Appl.  
\vol 17 
\yr 1972 
\pages 113--115
\endref 

\ref\n{Z} \by M. Zippin \paper Banach spaces with separable duals
 \jour Trans. Amer. Math. Soc. \vol 310 \yr 1988 \pages 371--379
 \endref


\endRefs
\enddocument